\documentclass[reqno]{amsart}

\usepackage{amsmath, amsthm, amscd, amsfonts, amssymb, graphicx, color}

\usepackage[bookmarksnumbered, colorlinks, plainpages]{hyperref}

\usepackage{lipsum}
\usepackage{fancyhdr}

\usepackage{algorithm,algorithmic}

\newtheorem{theorem}{Theorem}[section]
\newtheorem{lemma}[theorem]{Lemma}
\newtheorem{proposition}[theorem]{Proposition}
\newtheorem{corollary}[theorem]{Corollary}
\theoremstyle{definition}
\newtheorem{definition}[theorem]{Definition}
\newtheorem{example}[theorem]{Example}

\theoremstyle{remark}

\numberwithin{equation}{section}

\begin{document}
\title{More About Continuous GABOR FRAMES ON LOCALLY COMPACT ABELIAN GROUPS}
\author{Z. Hamidi}
\address{Department of Pure Mathematics, Azarbaijan shahid
 madani university, P.O.Box 53714-161, Tabriz, Iran.}
\email{zohre.hamidi@azaruniv.ac.ir}
\author{F. Arabyani-Neyshaburi}
\address{Department of Pure Mathematics, Ferdowsi University
of Mashhad, P.O.Box 1159-91775, Mashhad, Iran.}
\email{fahimeh.arabyani@gmail.com}
\author[R. A. Kamyabi-Gol]{R. A.  Kamyabi-Gol}
\address{Department of Mathematics, Faculty of Math, Ferdowsi University of Mashhad  and Center of Excellence in Analysis on Algebraic Structures (CEAAS), P.O.Box 1159-91775, Mashhad, Iran.}
\email{kamyabi@um.ac.ir}
\author{M. H. Sattari}
\address{Department of Pure Mathematics, Azarbaijan shahid
 madani university, P.O.Box 53714-161, Tabriz, Iran.}
\email{sattari@azaruniv.edu}



\subjclass[2000]{aa2bb; cc3dd.}

\keywords{LCA group, Zak transform, Continuous Gabor frame, Fibrization thechnique.}
\begin{abstract}
For a second countable locally compact abelian (LCA) group $G$, we study some necessary and sufficient conditions to generate continuous Gabor frames for $L^{2}(G)$.  To this end, we reformulate the generalized Zak transform proposed by Grochenig  in the case of integer-oversampled lattices, however our formulation  rely on the assumption that both translation and modulation groups are only closed subgroups. Moreover, we discuss the possibility of such generalization and apply several examples to demonestrate the necessity of standing conditions in the results. Finally, by using the generalized Zak transform and fiberization technique,  we obtain some characterization of continuous  Gabor frames for $L^{2}(G)$ in term of a family of frames in $l^{2}(\widehat{H^{\perp}})$ for a closed co-compact subgroup $H$ of $G$.
\end{abstract}

\maketitle

\section{Introduction}
The Zak transform is one of the  fundamental tools in both pure and applied mathematics,
that was originally introduced by Gelfand \cite{Gelfand} duo to some problems in differential equations.
This transform was studied  by weil   on locally compact abelian(LCA) groups \cite{W} and by Zak   in solid state physics \cite{Zak}. Later  it has been developed by many authors  for    identifying and characterizing of Gabor frames on $L^{2}(G)$ \cite{1, B-E, B-R, G}. The most research works in this regard are associated with a discrete, co-compact (uniform lattice) subgroup.
 In particular, Grochenig \cite{G}
presented some aspects of  Zak transform   to analyze  uniform lattice    Gabor frames. However, in recent years
this aspect  has been extended to closed subgroups  for  the characterization of continuous Gabor frames.  Indeed,  by considering a closed subgroup $H$ of an LCA group $G$ and applying the Zak transform associated to $H$   some equivalent conditions for the existence of continuous  Gabor frames as the form $\lbrace E_{\gamma} T_{\lambda} g \rbrace_{ \lambda \in H,\gamma \in H^{\perp}}$ have been obtained \cite{1, J-L}.

 The main purpose of this paper,  is to obtain some characterization results of   continuous  Gabor frames as   $\lbrace E_{\gamma} T_{\lambda} g \rbrace_{ \lambda \in \Lambda,\gamma \in \Gamma}$ on $L^{2}(G)$, for closed subgroups $\Lambda \subseteq G$  and $\Gamma \subseteq \widehat{G}$. To this end,  we  extend and reformulate
  the idea of  integer oversampling  for uniform lattices in \cite{G}. However,   our formulation  rely on the assumption that both translation and modulation groups are only closed subgroups and  remove some other  limited conditions.
    Moreover, we discuss the possibility of such generalization. Finally, by using the generalized Zak transform and fiberization method,  we give  some characterizations of continuous  Gabor frames for $L^{2}(G)$ in term of a family of frames in $l^{2}(\widehat{H^{\perp}})$ for a closed co-compact subgroup $H$ of $G$.

This paper is organized as follows. In Section $2$ we present  some   basic facts  of locally compact abelian groups and
the  required definitions of  continuous frame theory. Then we get a sufficient condition for the existence of  continuous Gabor frames in $ L^{2}(G)$. Section $3$ is devoted to investigate  some equivalent conditions for continuous Gabor frames in $ L^{2}(G)$ via generalized Zak transform.   Finally,  in section $4$, we use  the fiberization method to  survey   a relationship between  continuous  Gabor frames in $L^{2}(G)$ and a family of frames in $l^{2}(\widehat{H^{\perp}})$ for a closed co-compact subgroup $H$ of $G$.


\section{Notations and preliminaries}

Throughout this paper, let $G$ be a second countable locally compact abelian (LCA) group. It is known that such a group always carries a translation invariant regular Borel measure called a Haar measure, denoted by $\mu_{G}$, and is unique up to a positive constant.
We will use the addition as the group operation and  equip  discrete groups  with the counting measure.
Let $\widehat{G}$ denote the dual group of $G$, then
the famous duality theorem of Pontrjagin says that the character group of $\widehat{G}$ is topologically isomorphic with $G$, i.e., $\widehat{\widehat{G}} \cong G$.
The Fourier transform $\mathcal{F} :L^{1}(G)\longrightarrow C_{0} (\widehat{G})$,  is defined by
\begin{eqnarray*}
\mathcal{F} f (\xi) = \widehat{f} (\xi) = \int_{G} f(x)\overline{\xi} (x)   d\mu_{G} (x) \ \  \ (\xi \in \widehat{G}).
\end{eqnarray*}
 By the Fourier inversion, we can recover a function from its Fourier transform. Several different Fourier inversion theorem there exist. One of the most important states that if $f \in L^{1}(G)$ and $\widehat{f} \in L^{1}(\widehat{G})$, then
 \begin{eqnarray*}
 \mathcal{F}^{-1} \widehat{f}(x)= f(x)= \int_{\widehat{G}} \widehat{f}(\xi) \xi(x) d \mu_{\widehat{G}}(x)  \ \ \ ( a.e \ x \in G).
 \end{eqnarray*}
The Fourier transform can be extended from $L^{1}(G) \cap L^{2} (G)$ to an isometric isomorphism between $L^{2}(G) $ and $L^{2}(\widehat{G})$,  known as the Plancherel transform. See \cite{F, E-K-A, E-K}.

The mathematical theory for Gabor analysis in $L^{2}(G)$ is based on two classes of operators on $L^{2}(G)$. The translation  by $ \lambda \in G$, denoted by $T_{\lambda}$ and  is defined as $T_{\lambda} f (x)= f \left( x-\lambda\right)$, for all   $ x \in G$.
 Also, the modulation  by $ \gamma \in \widehat{G}$,  $E_{\gamma}$, defined by $E_{\gamma}f(x)=\gamma(x)f(x)$, for $x \in G$.
These classes of operators are unitary on $L^{2}(G)$ and satisfy the following relations;
 \begin{eqnarray*}
 T_{\lambda} E_{\gamma} = \overline{\gamma(\lambda)} E_{\gamma} T_{\lambda}, \ \ \mathcal{F}T_{\lambda}= E_{- \lambda}\mathcal{F} \  and \ \mathcal{F} E_{\gamma}= T_{\gamma} \mathcal{F}.
 \end{eqnarray*}
For a subset $\Lambda$ of an LCA group $G$, it's annihilator defined by
\begin{equation*}
\Lambda^{\perp}:= \lbrace \gamma \in \widehat{G}; \gamma (\lambda) =1 ,          \ for \ all \  \lambda \in \Lambda \rbrace.
\end{equation*}
 The annihilator is a closed subgroup of $\widehat{G}$ and if $\Lambda$ is a closed subgroup, then it is proved that $\widehat{\Lambda}\cong \frac{\widehat{G}}{\Lambda^{\perp}}$ and $\widehat{(\frac{G}{\Lambda})}\cong \Lambda^{\perp}$ \cite{F}. These relations show that for a closed subgroup $\Lambda$ the quotient $\frac{G}{\Lambda}$ is compact if and only if $\Lambda^{\perp}$ is discrete. See \cite{F, E-K-A, R} for more details.

 We also remind  the reader of Weil's formula; it relates integrable functions over $G$ with integrable functions on the quotient space $ \frac{G}{\Lambda}$ when $ \Lambda $ is a closed subgroup of $G$. For a closed subgroup $\Lambda$ of G we let $\pi_{\Lambda}: G\longrightarrow {\frac{G}{\Lambda}} ,  \pi_{\Lambda} (x) = x +\Lambda$ to be the canonical map from $G$ onto $\frac{G}{\Lambda}$. If $f \in L^{1}(G)$, then $\dot{x} :=\pi_{\Lambda} (x)$, defined almost everywhere on  $\frac{G}{\Lambda}$,  is integrable. Furthermore, when two of the Haar measures on $G$, $\Lambda$ and  $\frac{G}{\Lambda}$ are given, then the third can be normalized so that
\begin{eqnarray}\label{1}
 \int_{G} f(x) d \mu_{G}(x)= \int_{\frac{G}{\Lambda}} \int_{\Lambda} f(x+\lambda)d \mu_{\Lambda} (\lambda)d\mu_{\frac{G}{\Lambda}} (\dot{x}).
\end{eqnarray}
If  \eqref{1} holds, then the respective dual measures on $\widehat{G},  \Lambda^{\perp} \cong \widehat{\frac{G}{\Lambda}}$ and $ \frac{\widehat{G}}{\Lambda^{\perp}} \cong \widehat{\Lambda}$ satisfy
\begin{eqnarray}\label{2}
 \int_{\widehat{G}} \widehat{f}(\xi) d \mu_{\widehat{G}}(\xi) = \int_{\frac{\widehat{G}}{\Lambda^{\perp}}} \int_{\Lambda^{\perp}} \widehat{f}(\xi + \gamma)d \mu_{\Lambda^{\perp}} (\gamma)d\mu_{\widehat{G}/\Lambda^{\perp}} (\dot{\xi}).
\end{eqnarray}
Hence, if two of the measures on $G, \Lambda, {\frac{G}{\Lambda}}, \widehat{G}, \Lambda^{\perp}$ and $\frac{\widehat{G}}{\Lambda^{\perp}}$ are given, and these two are not dual measure, then by requiring Weil's formulas $\eqref{1}$ and $\eqref{2}$, all other measures are uniquely determined.
 For a closed subgroup  $\Lambda$ of $G$, a Borel section or a fundamental domain is a Borel measurable subset $X$ of $G$ such that every $y \in G$ can be uniquely written as $y= \lambda + x$, where $\lambda \in \Lambda $ and $ x \in X $. We always equip the Borel section $X$ of $G$ and the Haar measure $\mu_{G} \vert_{X}$. In \cite{B-R} it is shown that the mapping $x \longmapsto x+H$ from $(X, \mu_{G})$ to $(\frac{G}{H}, \mu_{\frac{G}{H}})$ is measure-preserving, and the mapping $Q(f) = f^{\prime}$ defined by
 \begin{equation}\label{measure}
 f^{\prime}(x+H)= f(x), \ \ \  x+H \in \frac{G}{H}, \ x \in X
 \end{equation}
 is an isometry from $ L^{2}(X , \mu_{G})$ onto $ L^{2}(\frac{G}{H}, \mu_{\frac{G}{H}})$ \cite{B-R}.

  Assume that $\Lambda$ is a discrete subgroup. It follows that $\mu_{G} (X)$ is finite if and only if, $\Lambda$ is co-compact, i.e., $ \Lambda$ is a uniform lattice  \cite{B}. For more information of harmonic analysis on locally compact abelian groups, we refer the reader to the classical books  \cite{ F, E-K-A, E-K, R}.

\subsection{Frame theory}
The central aspect  of this paper  is related to  continuous frames that was  introduced by Ali, Antoine and Cazeau \cite{A-C}. In the following we give the basic definitions and notations of continuous frames.
\begin{definition}
 Let $H$ be a complex Hilbert space, and let $(M, \sum_{M} , \mu_{M})$ be a measure space, where $\sum_{M}$ denotes the $\sigma$-algebra and $ \mu_{M}$ the non-negative measure. A family of vectors $\lbrace f_{k} \rbrace_{k \in M}$ is called a frame for $H$ with respect to $(M, \sum_{M} , \mu_{M})$ if

 (a) the mapping $M \longrightarrow \mathbb{C}$, $k\longmapsto \langle f, f_{k} \rangle$ is measurable for all $f \in \textit{H}$, and

 (b) there exist constants $A,B > 0$ such that
\begin{eqnarray}\label{3}
 A \parallel f \parallel^{2} \leq \int_{M} \vert \langle f  , f_{k} \rangle \vert^{2}  d \mu_{M} (k) \leq B\parallel f \parallel^{2}, \ \  \  (f \in \textit{H}).
\end{eqnarray}
 \end{definition}
 The constants A and B in $\eqref{3}$  are called frame bounds.
If $\lbrace f_{k} \rbrace_{k \in M}$ is weakly measurable and the upper bound in inequality  $\eqref{3}$ holds, then $\lbrace f_{k} \rbrace_{k \in M}$ is said to be a Bessel family with bound B. A frame $\lbrace f_{k} \rbrace_{k \in M}$ is said to be tight if we can choose $A=B$; if  furthermore $A=B=1$, then $\lbrace f_{k} \rbrace_{k \in M}$ is called a Parseval frame.
Also, a family $\lbrace f_{k} \rbrace_{k \in M}$ is called  minimal if $f_{j} \notin \overline{span} \lbrace f_{k} \rbrace_{k \neq j}$ for all $ j \in M$.

To a Bessel family $\lbrace f_{k} \rbrace_{k \in M}$ for $\textit{H}$, we associate the synthesis operator $T: L^{2}(M , \mu_{M}) \longrightarrow \textit{H}$ defined by
$T\lbrace c_{k} \rbrace_{k \in M} =\int_{M} c_{k} f_{k}  d \mu_{M} (k)$
 that  integral is defined in weak sense and is a bounded linear operator.

 Its adjoint operator $T^{*} : \textit{H} \longrightarrow L^{2}(M , \mu_{M})$  is called the $\textit{analysis operator}$, and  is given by
$T^{*}f = \lbrace \langle f , f_{k} \rangle \rbrace_{k \in M}$
The frame operator $S: \textit{H} \longrightarrow \textit{H}$ is defined as $S= TT^{*}$. We remark that the frame operator is the unique operator satisfying
\begin{equation*}
 \langle Sf , g \rangle = \int_{M} \langle f , f_{k} \rangle \langle f_{k} , g \rangle \ d \mu_{M}(k), \ \   \   (f,g \in H)
\end{equation*}
and   is a well-defined, bounded and self-adjoint for any Bessel system $\lbrace f_{k} \rbrace_{k \in M}$. Also,  it is invertible if and only if $\lbrace f_{k} \rbrace_{k \in M}$ is a frame \cite{A-C, R-N}.

Let $P$ be an index set, $g_{p}\in L^{2}(G)$ for all $p\in P$ and $H$ be a closed co-compact subgroup of $G$.
The translation invariant system generated by $\{g_{p}\}_{p\in P}$ with translation along the closed co-compact subgroup $H$ is as $\{T_{h}g_{p}\}_{h\in Hp\in P}$. Also, let for a topological space $T$, the Borel algebra of $T$ is denoted by $B_{T}$.
Then,    consider the following standing assumptions of \cite{J-L, J-S};

(I) $( P  , \sum_{P } , \mu_{P})$ is a $\sigma$-finite measure space,

(II) the mapping $p \longmapsto g_{p}, (P ,  \sum_{P}) \longrightarrow (L^{2}(G), B_{L^{2}(G)})$ is measurable,

(III) the mapping $(p , x ) \longmapsto  g_{p}(x) , (P \times G , \sum_{P} \bigotimes B_{G}) \longrightarrow ( \mathbb{C} , B_{\mathbb{C}})$  is measurable.

The family $\lbrace g_{p} \rbrace_{p \in P}$ is called admissible or, when $g_{p}$ is clear from the context, simply it is said that the measure space $P$ is admissible. The nature of these assumptions are discussed in  \cite{ J-S}. Observe that any closed subgroup $P_{j}$ of $G$ ( or $\widehat{G}$) with the Haar measure is admissible if $p \longmapsto g_{p}$ is continuous, e.g., if $g_{p} = T_{p}g$ for some function $g \in L^{2}(G)$.

A Gabor system in $L^{2}(G)$ with generator $ g \in L^{2}(G)$ is a family of functions of the form
$\lbrace E_{\gamma} T_{\lambda} g \rbrace_{ \lambda \in \Lambda,\gamma \in \Gamma}$,  where $\Gamma \leq \widehat{G}$ and  $\Lambda \leq G$.

If both   $\Lambda$ and $\Gamma$  are closed and co-compact subgroups, the family $ \lbrace E_{\gamma} T_{\lambda} g \rbrace_{\gamma \in \Gamma, \lambda \in \Lambda}$ is called a $\textit{co-compact Gabor system}$; if only one of the sets  $\Lambda$ or $\Gamma$ are closed and co-compact subgroup,  the Gabor system is called $\textit{ semi co-compact} $ and
 it is said to be a uniform lattice Gabor system whenever  both   $\Lambda$  and $\Gamma$ are uniform lattices.

In the following, we are going to derive some sufficient  condition for $\lbrace E_{\gamma} T_{\lambda}g \rbrace_{\lambda \in \Lambda, \gamma \in \Gamma}$ being a frame for $L^{2}(G)$. First we need the following  result.

\begin{lemma} \cite{J-L}\label{JacobL}
Let $g \in L^{2}(G)$ and $\Gamma \subseteq \widehat{G}$ be a closed subgroup. Then
\begin{eqnarray*}
&& \int_{\Gamma} \vert\left\langle f, E_{\gamma} T_{\lambda} g\right\rangle \vert^{2} d \mu_{\Gamma}(\gamma) \\&=& \int_{G} \int_{\Gamma^{\perp}} f(x) \overline{f(x- \alpha)} \overline{ T_{\lambda} g(x)} T_{\lambda}g(x- \alpha) d \mu_{\Gamma^{\perp}} (\alpha) d \mu_{G}(x)
\end{eqnarray*}
for any $f \in C_{c} (G)$.
\end{lemma}
\begin{proposition}
Let 
$ (\Lambda , \mu_{\Lambda}) \subseteq G$ be an admissible measure space and $\Gamma \subseteq \widehat{G}$ be a closed and co-compact subgroup. Also, let for all $\alpha \in \Gamma^{\perp}$ we have $suppg \bigcap supp T_{\alpha} g = \emptyset$,  up to a set of measure zero in $G$. If there exist constants $A, B > 0$ such that
\begin{eqnarray}\label{cndi}
A \leq \int_{\Lambda} \vert T_{\lambda} g(x) \vert ^{2} d \mu_{\Lambda} (\lambda) \leq B,     \ \ \ a.e. \ (x \in G).
\end{eqnarray}
Then $\lbrace E_{\gamma} T_{\lambda}g \rbrace_{\lambda \in \Lambda, \gamma \in \Gamma}$ is a frame for $L^{2}(G)$.
\end{proposition}
\begin{proof}
Since $\Gamma$ is closed and  co-compact  so $\Gamma^{\perp}$ is a   discrete subgroup of $G$. Thus by the assumption and Lemma \ref{JacobL}  we have
\begin{eqnarray*}
&&\int_{\Gamma} \vert \langle f, E_{\gamma} T_{\lambda}g \rangle \vert^{2} d \mu_{\Gamma} (\gamma) \\
& =& \int_{G} \sum_{\alpha \in \Gamma^{\perp}} f(x) \overline{f(x- \alpha)} \overline{ T_{\lambda}g(x)} T_{\lambda} g(x-\alpha) d \mu_{G}(x)\\
& =& \int_{G} \vert f(x) \vert^{2} \vert T_{\lambda}g(x) \vert^{2} d \mu_{G}(x) \\
&+& \int_{G} \sum_{ \alpha \neq 0 , \alpha \in \Gamma^{\perp}} f(x) \overline{f(x- \alpha)} \overline{ T_{\lambda}g(x)} T_{\lambda} g(x-\alpha) d \mu_{G} (x) \\
 & =& \int_{G} \vert f(x) \vert^{2} \vert T_{\lambda}g(x) \vert^{2} d \mu_{G}(x),
\end{eqnarray*}
for every $f \in C_{c}(G)$.
Therefore,
\begin{eqnarray*}
\int_{\Lambda} \int_{\Gamma} \vert \langle f, E_{\gamma} T_{\lambda}g \rangle \vert^{2} d \mu_{\Gamma} (\gamma) d \mu_{\Lambda}(\lambda) = \int_{G} \vert f(x) \vert^{2} \int_{\Lambda} \vert T_{\lambda}g(x) \vert^{2} d \mu_{\Lambda}(\lambda) d \mu_{G}(x)
\end{eqnarray*}
and consequently by using $(\ref{cndi})$
\begin{eqnarray}\label{4}
 A\parallel f \parallel^{2} \leq \int_{\Lambda} \int_{\Gamma} \vert \langle f, E_{\gamma} T_{\lambda}g \rangle \vert^{2} d \mu_{\Lambda}(\lambda) \ d \mu_{\Gamma}(\gamma) \leq B \parallel f \parallel^{2}
\end{eqnarray}
for all $f \in C_{c} (G)$. In addition, $\Lambda$ and $\Gamma$ are $\sigma$-finite measure spaces and  $ \overline{C_{c}(G)} = L^{2}(G)$. Thus,   Proposition 2.5 of \cite{R-N} implies that  the inequality \eqref{4} holds for all $f \in L^{2}(G)$, i.e., $\lbrace E_{\gamma} T_{\lambda}g \rbrace_{\lambda \in \Lambda, \gamma \in \Gamma}$ is a frame for  $L^{2}(G)$.
\end{proof}
\section{Continuous Zak transform and construction of Gabor frames}
In this section, we deal with the generalized Zak transform related to continuous Gabor systems on LCA groups.  We  extend  the idea of  integer oversampling to closed subgroups and remove some limited assumptions. Moreover, we  discuss the existence conditions and provide equivalent conditions for a continuous Gabor system to be a frame family, orthonormal basis, complete or a  minimal family.
\begin{definition}
Let $G$ be an LCA group and $\Lambda$ be a closed subgroup of $G$. For each $f \in C_{c}(G)$, the mapping $Z_{\Lambda}f$, defined on $G \times \widehat{G}$  as
\begin{equation*}
Z_{\Lambda} f(x,\xi) = \int_{\Lambda} f(x +\lambda) \xi(\lambda) d m_{\Lambda}(\lambda)
\end{equation*}
 is called the continuous Zak transform of $f$.
\end{definition}
It is known that, the continuous Zak transform can be extended to a unitary operator from $ L^{2}(G)$ onto $ L^{2} (M_{\Lambda}) $, where $M_{\Lambda} := \frac{G}{\Lambda} \times  \frac{\widehat{G}}{\Lambda^{\perp}}$.
The next lemma states the basic properties of continuous Zak transform.
\begin{lemma}\cite{1}\label{lem1.31}
Let  $\Lambda$ be a closed subgroup of $G$ and $f\in L^{2}(G)$. Then
\begin{itemize}
\item[(i)]Quasi-periodicity: for all $a \in G$ and $\gamma \in \Lambda^{\perp}$
\begin{equation*}
Z_{\Lambda}f\left(x+ a, \xi\right) = \overline{\xi(a)} Z_{\Lambda}f\left( x, \xi\right) , \ \ \   Z_{\Lambda}f\left( x, \xi + \gamma\right)=Z_{\Lambda}f\left( x, \xi\right).
\end{equation*}
\item[(ii)] Diagonalization: If $\left( \gamma, \lambda \right) \in \left( \Lambda^{\perp} \times \Lambda \right)$, then $E_{\gamma} T_{\lambda} f \in L^{2}(G)$ and
$Z_{\Lambda}E_{\gamma} T_{\lambda}f= E_{\lambda, \gamma}Z_{\Lambda}f$,
where $E_{\lambda, \gamma} (x , \omega)= \gamma(x) \omega(\lambda)$ for all $(x, \omega) \in G \times \widehat{G}$.
\end{itemize}
\end{lemma}

Now, let $G$ be an LCA group,  $\Lambda \subseteq G $ and $\Gamma\subseteq \widehat{G}$ be closed subgroups.
Also, let there exists a closed subgroup $H \leq \Lambda$ so that $ H^{\perp} \leq \Gamma $ and $ \frac{\Lambda}{H} $ and $ \frac{\Gamma}{H^{\perp}}$ are countable (or finite).
In this case we can choose $\lambda_{i} \in \Lambda$ so that $\Lambda = \cup_{i=1}^{\infty} \left( \lambda_{i} + H\right)$ and each coset of $\frac{\Lambda}{H}$ contains only one $\lambda_{i}$. Moreover, there exist $\gamma_{j}\in \Gamma$,  so that $\Gamma= \cup_{j=1}^{\infty} \left( \gamma_{j}+  H^{\perp}\right) $ and each coset of $\frac{\Gamma}{H^{\perp}}$ contains only one $\gamma_{j}$. Then the frame operator of  the Gabor system $\lbrace E_{\gamma} T_{\lambda} g \rbrace_{\gamma \in \Gamma, \lambda \in \Lambda}$, for a well-fitted window   function  $g$ on $G$, can be written as follows
\begin{eqnarray*}
Sf &=& \int_{\Gamma}\int_{\Lambda}  \langle f , E_{\gamma} T_{\lambda} g \rangle E_{\gamma} T_{\lambda}g
d\mu_{\Lambda}(\lambda) d\mu_{\Gamma}(\gamma)\\
&=& \sum_{i=1}^{\infty} \sum_{j=1}^{\infty}  \int_{H^{\perp}}\int_{H} \langle f , E_{\omega} T_{h} T_{\lambda_{i}} E_{\gamma_{j}} g \rangle E_{\omega} T_{h} T_{\lambda_{i}} E_{\gamma_{j}} g d \mu_{H}(h) d \mu_{H^{\perp}}(\omega)\\
&=&  \sum_{i=1}^{\infty} \sum_{j=1}^{\infty} \int_{H^{\perp}}\int_{H}  \langle f , E_{\omega} T_{h} g_{ij} \rangle E_{\omega} T_{h} g_{ij} d \mu_{H}(h) d \mu_{H^{\perp}}(\omega)
\end{eqnarray*}
for all $f \in L^{2}(G)$ where
\begin{equation}\label{gij}
g_{ij} = T_{\lambda_{i}} E_{\gamma_{j}} g.
\end{equation}
Thus
\begin{eqnarray*}
Z_{H}Sf &=& \sum_{i=1}^{\infty} \sum_{j=1}^{\infty}  \int_{H^{\perp}} \int_{H}\langle Z_{H}f , Z_{H} \left(  E_{\omega} T_{h} g_{ij} \right)  \rangle Z_{H} \left(  E_{\omega} T_{h} g_{ij} \right)  d \mu_{H}(h) d \mu_{H^{\perp}}(\omega) \\
&=&\sum_{i=1}^{\infty} \sum_{j=1}^{\infty}  \int_{H^{\perp}}\int_{H} \langle Z_{H}f , E_{\omega , h} Z_{H} g_{ij} \rangle E_{\omega , h} Z_{H} g_{ij}   d \mu_{H}(h) d \mu_{H^{\perp}}(\omega)  \\
&=& \sum_{i=1}^{\infty} \sum_{j=1}^{\infty} \int_{H^{\perp}} \int_{H} \widehat{\left( Z_{H}f . \overline{Z_{H} g_{ij}} \right) } (h,\omega) E_{\omega , h} Z_{H} g_{ij}   d \mu_{H}(h) d \mu_{H^{\perp}}(\omega)   \\
&=& \sum_{i=1}^{\infty} \sum_{j=1}^{\infty} Z_{H}f . \overline{Z_{H} g_{ij}} . Z_{H} g_{ij}\\
&=& \left( \sum_{i=1}^{\infty} \sum_{j=1}^{\infty}  \mid Z_{H} g_{ij} \mid^{2} \right) Z_{H}f.
\end{eqnarray*}
The forthcoming theorem, which collects the above computations, shows that   the Zak transform on $H$ diagonalize the Gabor frame operator of $\lbrace E_{\gamma} T_{\lambda} g \rbrace_{\lambda \in \Lambda, \gamma \in \Gamma}$. Moreover,  the spectrum of the frame operator   equals the range of  $\sum_{i=1}^{\infty} \sum_{j=1}^{\infty}  \mid Z_{H} g_{ij} \mid^{2}$.
\begin{theorem}\label{diagonalize}
Let  $g$,  $\Lambda$, $\Gamma$ and $S$  be as above and 
 there exists a closed subgroup $H$ of $G$ so that
\begin{eqnarray}\label{existence H}
H \leq \Lambda,   \textit{ and }  H^{\perp} \leq \Gamma.
\end{eqnarray}
Moreover, assume that    $\frac{\Lambda}{H}$  and    $\frac{\Gamma}{H^{\perp}}$  are countable.
Then, we obtain  $Z_{H}SZ_{H}^{-1}F = \left( \sum_{i=1}^{\infty} \sum_{j=1}^{\infty}  \mid Z_{H} g_{ij} \mid^{2} \right)F$, for all $F \in L^{2} (M_{H})$.
\end{theorem}
As a special case of Theorem \ref{diagonalize} we record the following corollaries.

\begin{corollary}\label{3.4..}
Let $G$ be an LCA group, $g\in L^{2}(G)$,  $H, \Lambda \leq G $ and $\Gamma \leq \widehat{G}$ be closed subgroups.   Then,
\begin{itemize}
\item[(i)]
$Z_{H}\left( E_{\gamma} T_{\lambda} g\right) (x, \omega)= \gamma( \lambda) E_{\lambda , \gamma}(x, \omega) Z_{H} g(x, \gamma +\omega)$, for all $\lambda\in \Lambda$, $\gamma\in \Gamma$  and a.e.   $ (x , \omega) \in G \times \widehat{G} $.
\item[(ii)] If  the  closed subgroup $H$ of $G$ satisfies $(\ref{existence H})$, then $Z_{H}\left( E_{\gamma} T_{\lambda} g\right)=  E_{\lambda , \gamma} Z_{H} g$, for all $\lambda\in \Gamma^{\perp}$ and  $\gamma\in \Lambda^{\perp}$.
\end{itemize}
\end{corollary}
\begin{proof}
 To show $(i)$, suppose that $g\in L^{2}(G)$ then
\begin{eqnarray*}
&& Z_{H}\left( E_{\gamma} T_{\lambda} g \right) (x, \omega) = \int_{H} E_{\gamma} T_{\lambda} g (x+h) \omega(h) d \mu_{H} (h)  \\
&=& \int_{H} g(x+h- \lambda) \gamma(x) \gamma(h) \omega(h) d \mu_{H} (h)   \\
&=& \int_{H} g(x+h) (\gamma +\omega)(h) (\gamma+ \omega) (\lambda) \gamma(x) d \mu_{H} (h)  \\
&=& (\gamma +\omega) (\lambda) \gamma(x) Z_{H} g(x, \gamma+ \omega) \\
& =& \gamma( \lambda) E_{\lambda , \gamma}(x, \omega) Z_{H} g(x, \gamma+ \omega),
\end{eqnarray*}
for all  $\lambda\in \Lambda$, $\gamma\in \Gamma$, $x\in G$ and $\omega\in \widehat{G}$. The proof of $(ii)$ is similar, only  it is sufficient to  note that in the above computations we have   $\gamma( \lambda) =\gamma(h)=1$, for all $h\in H, \lambda\in \Gamma^{\perp}, \gamma\in \Lambda^{\perp}$, by the assumption.
\end{proof}
\begin{corollary}
Let $G$,  $g$, $\Lambda$,  $\Gamma$, and $H$  satisfy in conditions of Theorem \ref{diagonalize},
and take the  sequence  $ \lbrace g_{ij} \rbrace_{ i  , j =1}^{\infty}$ as in $(\ref{gij})$. Then,  the following statements hold:
\begin{itemize}
\item[(i)] $ \lbrace E_{\gamma} T_{\lambda} g \rbrace_{\lambda \in \Lambda , \gamma \in \Gamma}$  is a frame  for $ L^{2} (G)$  if and only if there exist constants  $A$, $B$ so that $ A \leq \left( \sum_{i=1}^{\infty} \sum_{j=1}^{\infty}  \mid Z_{H} g_{ij}(x , \omega) \mid^{2} \right) \leq B $, for a.e. $ (x , \omega) \in G \times \widehat{G}$.
\item[(ii)] $ \lbrace E_{\gamma} T_{\lambda} g \rbrace_{\lambda \in \Lambda , \gamma \in \Gamma}$ is a Parseval frame  for $ L^{2} (G) $  if and only if
$\sum_{i=1}^{\infty} \sum_{j=1}^{\infty}  \mid Z_{H} g_{ij}(x , \omega) \mid^{2} =1 $, for a.e. $ (x , \omega) \in G \times \widehat{G}$.
\item[(iii)] $ \lbrace E_{\gamma} T_{\lambda} g \rbrace_{\lambda \in \Lambda , \gamma \in \Gamma}$ is an orthonormal basis     for $ L^{2} (G) $  if and only if
$ \Vert g \Vert_{L^{2}(G)} =1$ and $ \sum_{i=1}^{\infty} \sum_{j=1}^{\infty}  \mid Z_{H} g_{ij}(x , \omega) \mid^{2} =1$, for a.e. $ (x , \omega) \in G \times \widehat{G}$.
\end{itemize}
\end{corollary}

\begin{proposition}\label{3.6prop./}
Let $G$ be an LCA group, $g \in L^{2}(G)$, and $H$  be a closed subgroup of $G$. Then,  the following statements hold:
\begin{itemize}
\item[(i)]   $ \lbrace E_{\gamma} T_{\lambda} g \rbrace_{\lambda \in H , \gamma \in H^{\perp}}$  is a complete system  in $ L^{2} (G)$ if and only if   $ Z_{H}g \neq 0 $, a.e..
\item[(ii)] If $H$ is a uniform lattice, then   $ \lbrace E_{\gamma} T_{\lambda} g \rbrace_{\lambda \in H , \gamma \in H^{\perp}}$ is a  minimal system  in  $L^{2} (G)$  if and only if $ \dfrac{1}{\overline{Z_{H}g}} \in L^{2}(M_{H})$.
\end{itemize}
\end{proposition}
\begin{proof}
$(i)$; Let $Z_{H}g \neq 0$ a.e.,
 to show the Gabor system $\lbrace E_{\gamma}T_{\lambda}g \rbrace_{\lambda \in H , \gamma \in H^{\perp}}$ is complete in $L^{2} (G)$,
it is sufficient to prove   $ \lbrace   E_{\lambda, \gamma}Z_{H} g \rbrace_{\lambda \in H , \gamma \in H^{\perp}}$ is complete in $ L^{2}(M_{H})$, by  the unitarity of $ Z_{H}$ and Lemma \ref{lem1.31} $(ii)$. To this end, let $\Phi \in L^{2}(M_{H})$ such that $\langle \Phi  ,  E_{\lambda, \gamma}Z_{H} g \rangle_{ L^{2}(M_{H})} = 0$, for all $\lambda \in H , \gamma \in H^{\perp}$.  So we can write
\begin{eqnarray*}
&& \int_{\frac{G}{H}} \int_{\frac{\widehat{G}}{H^{\perp}}} \Phi  ( \alpha , \beta) \overline{Z_{H}g(\alpha , \beta)}  \overline{ E_{\lambda, \gamma} ( \alpha , \beta) } \ d\mu_{\frac{G}{H}} (\dot{\alpha}) \ d \mu_{\frac{\widehat{G}}{H^{\perp}}} ( \dot{\beta})\\
&=& \langle   \Phi , E_{\lambda, \gamma} Z_{H}g \rangle_{ L^{2}(M_{H})} =0,
\end{eqnarray*}
for all $ \lambda \in H , \gamma \in H^{\perp} $.
Since $ \phi . \overline{Z_{H}g} \in L^{1}(M_{H})$ and  the functions in $ L^{1}(M_{H}) $ are uniquely determined by their Fourier coefficients. Thus,  we have $ \Phi . \overline{Z_{H}g}=0$, a.e..
On the other hand, by the assumption   $Z_{H}g\neq 0$ a.e. that means
$\Phi =0 $ a.e.,  proving the claim.

For the converse, suppose the Gabor system  $ \lbrace E_{\gamma} T_{\lambda} g \rbrace_{\lambda \in H , \gamma \in H^{\perp}}$  is  a complete family   in $ L^{2}(G) $. Then $\{E_{\lambda , \gamma} Z_{H} g\}_{\lambda \in H , \gamma \in H^{\perp}}$ is also   complete    in $ L^{2}(M_{H})$. To the contrary, let $ \Delta_{g} = \lbrace ( \alpha , \beta)  \in M_{H} : Z_{H}g( \alpha , \beta)  =0 \rbrace $ has positive measure. Put  $\Phi= \chi_{\Delta_{g}}$, that $ \chi_{\Delta_{g}}$ is the characteristic function on $ \Delta_{g}$. Then we obtain $ \langle \Phi , E_{\lambda , \gamma} Z_{H} g \rangle_{ L^{2}(M_{H})}=0$,
 for  all $\lambda \in H , \gamma \in H^{\perp}$, a contradiction.
So,  $ Z_{H}g\neq 0 $ a.e..


To show $(ii)$, first suppose that $ \dfrac{1}{\overline{Z_{H}g}} \in L^{2}(M_{H})$.
Since  $Z_{H}$ is surjective, there exists a   function  $ h \in L^{2} (G)$ such that $ \dfrac{1}{\overline{Z_{H}g}} = Z_{H}h$.
Hence, we have
\begin{eqnarray*}
&& \langle E_{\gamma} T_{\lambda} g , E_{\gamma^{\prime}} T_{\lambda^{\prime}} h \rangle_{ L^{2} (G)}\\
&=&  \langle    E_{\lambda , \gamma} Z_{H} g  ,   E_{\lambda^{'} , \gamma^{'}}Z_{H}h \rangle_{ L^{2} (M_{H})}\\
 &=& \langle    E_{\lambda , \gamma} Z_{H} g  ,   E_{\lambda^{'} , \gamma^{'}} \dfrac{1}{\overline{Z_{H}g}} \rangle_{ L^{2} (M_{H})} \\
&=&   \langle    E_{\lambda , \gamma},    E_{\lambda^{'} , \gamma^{'}} \rangle_{ L^{2} (M_{H})},
\end{eqnarray*}
for all $\lambda, \lambda^{\prime}\in H$ and $\gamma, \gamma^{\prime} \in H^{\perp}$. The assumption that $H$ is a uniform lattice implies  that $M_{H}$ is a compact group and hence
$\lbrace E_{\lambda , \gamma} \rbrace_{\lambda \in H, \gamma \in H^{\perp}}$ is an orthonormal basis for $L^{2}(M_{H})$.
Therefore, it follows from the above computations that  $\lbrace E_{\gamma} T_{\lambda} g \rbrace_{\lambda \in H , \gamma \in H^{\perp}}$ has a biorthogonal system and consequently it  is  minimal.
Conversely, let $ \lbrace E_{\gamma} T_{\lambda} g \rbrace_{\lambda \in H , \gamma \in H^{\perp}} $ is minimal so $ \lbrace E_{\lambda , \gamma} Z_{H}g \rbrace_{ \lambda \in H , \gamma \in H^{\perp}}$ is minimal and has a biorthogonal system as $ \lbrace \mathcal{\psi}_{\lambda , \gamma} \rbrace_{ \lambda \in H , \gamma \in H^{\perp}} \subseteq L^{2}(M_{H})$. Fixed $ \lambda_{0} \in H,  \gamma_{0} \in H^{\perp}$, then $ \mathcal{\psi}_{\lambda_{0},\gamma_{0}}. \overline{Z_{H}g} \in L^{1}(M_{H})$ and we have
\begin{eqnarray*}
&& \int_{\frac{G}{H}} \int_{\frac{\widehat{G}}{H^{\perp}}} \Phi  ( \alpha , \beta) \overline{Z_{H}g(\alpha , \beta)}  \overline{ E_{\lambda, \gamma} ( \alpha , \beta) } \ d\mu_{\frac{G}{H}} (\dot{\alpha}) \ d \mu_{\frac{\widehat{G}}{H^{\perp}}} ( \dot{\beta}) \\
&=& \langle \mathsf{\psi}_{\lambda_{0},\gamma_{0}} , Z_{H}g  E_{\lambda , \gamma} \rangle_{L^{2}(M_{H})} = \delta_{\lambda , \lambda_{0}} \delta_{\gamma , \gamma_{0}}
\end{eqnarray*}
on the other hand $ \langle E_{\lambda , \gamma} , E_{\lambda_{0} , \gamma_{0}} \rangle = \delta_{\lambda , \lambda_{0}} \delta_{\gamma , \gamma_{0}}$ for all $ \lambda \in H , \gamma \in H^{\perp}$. so $\mathrm{\psi}_{\lambda_{0},\gamma_{0}} \overline{Z_{H}g} = E_{\lambda_{0} , \gamma_{0}} \neq 0$ thus $ \mathtt{\psi}_{\lambda_{0},\gamma_{0}} = \dfrac{E_{\lambda_{0} , \gamma_{0}}}{\overline{Z_{H}g}} $. In particular, we obtain $ \mathbb{\psi}_{e,1} = \dfrac{1}{\overline{Z_{H}g}} \in L^{2}(M_{H})$.
\end{proof}
It is worth noticing that for closed subgroups $ \Lambda$,$\Gamma$ and $H$ which satisfy (\ref{existence H}), we have
$\Gamma^{\perp} \times \Lambda^{\perp} \subseteq H \times H^{\perp}\subseteq \Lambda \times \Gamma$. 
 So, if $Z_{H} \neq 0$ a.e. then the Gabor syatem $ \lbrace E_{\gamma} T_{\lambda} g \rbrace_{ \lambda \in \Lambda, \gamma \in \Gamma}$ is complete in $L^{2}(G)$.  However, as the following examples show, the Gabor system $ \lbrace E_{\gamma} T_{\lambda} g \rbrace_{ \lambda \in \Gamma^{\perp}, \gamma \in \Lambda^{\perp}}$ is not necessarily complete and in case $ \dfrac{1}{\overline{Z_{H}g}} \in L^{2}(M_{H})$, the Gabor system $\lbrace E_{\gamma} T_{\lambda} \phi \rbrace_{\lambda \in \Lambda , \gamma \in \Gamma}$ is not minimal, in general.
\begin{example}
Fix $ 0 < \alpha < 1$, set $g(x)= \vert x \vert^{\alpha} $ for $x \in \left[ - \frac{1}{2} , \frac{1}{2} \right] $. It is known  that the system $ \lbrace T_{n} E_{m} g \rbrace_{n,m \in \mathbb{Z}}$ is a Schauder basis for $L^{2}(G)$ (but not Riesz basis for $L^{2}(G)$), \cite{D-H}. So this system is minimal and complete.  Take $ \Lambda = \frac{1}{2} \mathbb{Z}$ , $\Gamma = \frac{1}{4} \mathbb{Z}$ and $H= \mathbb{Z}$, then the closed subgroups $H , \Lambda$ and $ \Gamma$ satisfy (\ref{existence H}). In addition, the Gabor system $ \lbrace E_{\gamma} T_{\lambda} \phi \rbrace_{ \lambda \in H, \gamma \in H^{\perp}}$ is complete and minimal. However, the Gabor system $ \lbrace E_{\gamma} T_{\lambda} \phi \rbrace_{\lambda \in \Lambda , \gamma \in \Gamma}$ is not minimal and  the Gabor system $\lbrace E_{\gamma} T_{\lambda} \phi \rbrace_{\gamma \in \Lambda^{\perp} , \lambda \in \Gamma^{\perp}}$ is not complete.
\end{example}
As another example, suppose that  $G= \mathbb{R}$ and consider the $Gaussian$ function $\phi (x)= e^{-\pi x^{2}}$. It was already proved in \cite{G-R} that  the Gabor system $ \lbrace E_{m \alpha} T_{n \beta} \phi \rbrace_{ m,n \in \mathbb{Z}}$ is complete for $ \alpha \beta \leq 1$ and incomplete for $ \alpha \beta > 1$. Moreover,  $ Z_{\alpha \mathbb{Z}} \phi \neq 0$ a.e. on $ \left[ 0 , \alpha \right) \times \left[ 0 , \frac{1}{\alpha} \right)$, for all non-zero $ \alpha \in \mathbb{R}$ \cite{H1}. Let $ H= 4 \mathbb{Z}$,  $\Lambda = 2 \mathbb{Z}$ and $ \Gamma = \frac{1}{8} \mathbb{Z}$, then  the closed subgroups $H$, $\Lambda$ and $\Gamma$ satisfy (\ref{existence H}). So, the Gabor system $ \lbrace E_{\gamma} T_{\lambda} \phi \rbrace_{ \lambda \in H, \gamma \in H^{\perp}}$ is complete, and so $Z_{H}\phi\neq 0$, a.e., by Proposition \ref{3.6prop./}. Although,
the Gabor system $ \lbrace E_{\gamma} T_{\lambda} \phi \rbrace_{\gamma \in \Lambda^{\perp} , \lambda \in \Gamma^{\perp}}$ is incomplete.


\subsection{The existence conditions}
In what follows, for two given closed subgroups $\Lambda$ and $\Gamma$, we discuss the existence of a closed subgroup $H$ which satisfies $(\ref{existence H})$ so that the quotient groups $ \frac{\Lambda}{H}$ and $\frac{\Gamma}{H^{\perp}}$ are countable.
We first  note that,  for every  LCA group $G$, the condition $ \Gamma^{\perp} \leq \Lambda$ is  necessary  for the existence of $H$.
Moreover obviously, for countable groups it can be considered as a necessary and sufficient condition.
Specially, for finite group $G = \mathbb{Z}_{L} = \lbrace 0,..., L-1 \rbrace$,  $L \in \mathbb{N}$, it is known that
 a closed subgroup in $G$  is as $\Lambda= N \mathbb{Z}_{\frac{L}{N}}$ where $N \in \mathbb{N}$ so that $N$ is a divisor of $L$. 
Also, consider closed subgroups $\Gamma = M \mathbb{Z}_{\frac{L}{M}}$ and $ H= R \mathbb{Z}_{\frac{L}{R}}$ so that  $M, R \in \mathbb{N}$, are some divisors of $L$. Then the condition (\ref{existence H}) is equivalent to
\begin{eqnarray}\label{5}
N \mid R \mid L \ \  and \  \  M \mid L/R.
\end{eqnarray}
More precisely, (\ref{5}) is the necessary and sufficient condition for the subgroup $H= R \mathbb{Z}_{\frac{L}{R}}$ of $G$ to satisfy (\ref{existence H}).

\begin{lemma}\label{3.7.lem}
Assume that $H \leq \Lambda$ are closed subgroups of $G$ so that $\frac{\Lambda}{H}$ is finite. Then $ \frac{\Lambda}{H} \cong \frac{ H^{\perp}}{\Lambda^{\perp}}$.
\end{lemma}
\begin{proof}
Applying Proposition 4.2.24 of \cite{R} and the fact that any finite group  is self-dual  we obtain
$\frac{\Lambda}{H} \cong \widehat{\left(\frac{\Lambda}{H}\right)}  \cong \frac{ H^{\perp}}{\Lambda^{\perp}}$.
\end{proof}
\begin{example}\label{ex3.5.}
Suppose that $G = Q_{p}$,  the $P$-adic numbers group, it is known that every non-trivial closed subgroup $H$ of $Q_{p}$, is open and compact. Hence, $ \frac{G}{H}$ is infinite and discrete, so for  every  non-trivial closed subgroups $H\leq \Lambda$ of $Q_{p}$ we imply that $\frac{\Lambda}{H}$ is both discrete and compact and consequently is finite, similarly $ \frac{\Gamma}{H^{\perp}}$  is finite as well. That means all closed subgroups $H$ with the property $\Gamma^{\perp} \leq H \leq \Lambda$, fulfill the condition  (\ref{existence H}) a long with finite quotient groups.
\end{example}
\begin{example}
Consider   $G= \mathbb{R} \times Z_{p}$, where $Z_{p}$ is the group of $p$-adic integers, and let $ \Lambda , \Gamma $ be two non-trivial closed subgroups of $G$ and $ \widehat{G}$, respectively. So $ \Lambda = \alpha \mathbb{Z} \times \Lambda_{2}$, for some  $ \alpha \in \mathbb{R}$ and $ \Lambda_{2}$ is a closed subgroup of $Z_{p}$. We show that  for   every closed subgroup $H$ which satisfies 
$\Gamma^{\perp} \leq H \leq \Lambda$,  the quotients $ \frac{\Lambda}{H}$ and $ \frac{\Gamma}{H^{\perp}}$ are finite. Indeed $H \leq \Lambda \lneq \mathbb{R} \times Z_{p}$ implies that $H = \alpha m \mathbb{Z} \times H_{2}$, for some $m \in \mathbb{Z}$ and $H_{2} \leq \Lambda_{2} \leq Z_{p}$. Hence a similar discussion as in Example \ref{ex3.5.} assures that $ \frac{\Lambda_{2}}{H_{2}}$ is finite and consequently
\begin{equation*}
\frac{\Lambda}{H} = \frac{\alpha \mathbb{Z}}{\alpha m \mathbb{Z}} \times \frac{\Lambda_{2}}{H_{2}}
\end{equation*}
is a finite group. Moreover, $ \Gamma^{\perp} \leq H \leq Z_{p}$ also follows that $ \frac{H}{\Gamma^{\perp}}$ is finite and so by Lemma \ref{3.7.lem} we obtain $ \frac{\Gamma}{H^{\perp}} \cong \frac{H}{\Gamma^{\perp}}$, i.e., $ \frac{\Gamma}{H^{\perp}}$ is finite, as well.
\end{example}

In the sequel,  we  show that, $ \Gamma^{\perp} \leq \Lambda$  is not a sufficient condition for the existence of desired $H$ with countable quotient groups, in general.
\begin{example}\label{ex3.11}
Let $ G= \mathbb{R}^{n}$ and $ \Lambda = \Gamma = \mathbb{R}  \times \mathbb{Z}^{n-1}$. Then  $\Gamma^{\perp}  \leq \Lambda$ and
 for every closed subgroup $H$ so that  $\Gamma^{\perp}  \leq H \leq  \Lambda$  we can write $H = H_{1} \times H_{2}$ where $H_{1} \leq \mathbb{R}$ and $H_{2} = \mathbb{Z}^{n-1}$. If $H_{1} \neq \mathbb{R}$, then $H_{1} = \alpha \mathbb{Z}$ for some $ \alpha \in \mathbb{R}$. Thus $ \frac{\Lambda}{H}$ and $ \frac{\Gamma}{H^{\perp}}$  are  uncountable. Moreover, if $H_{1}= \mathbb{R}$ i.e,. $H = \mathbb{R} \times \mathbb{Z}^{n-1}$, then $ \dfrac{\Gamma}{H^{\perp}}$
 is uncountable.
\end{example}

 Suppose $G$ is a compactly generated group  of Lie type, that is isomorphic to one of the form $ G= \mathbb{R}^n \times \mathbb{Z}^m \times \mathbb{T}^r \times F$ for $ n,m,r \in \mathbb{N}$ and  a finite abelian group $F$. Consider  closed subgroups $ \Lambda = \Lambda_{1} \times \Lambda_{2} \times \Lambda_{3} \times \Lambda_{4} \leq G$ and  $ \Gamma = \Gamma_{1} \times  \Gamma_{2} \times  \Gamma_{3} \times  \Gamma_{4} \leq \widehat{G}$ so that $ \Lambda_{1} = \Gamma_{1} = \mathbb{R} \times \mathbb{Z}^{n-1}$, and  $\Gamma_{i}=\Lambda_{i}^{\perp}$, for $2\leq i \leq 4$. Then for any closed subgroup $H$ so that  $\Gamma^{\perp}\leq H \leq \Lambda$ we have that $ H= H_{1} \times H_{2} \times H_{3} \times H_{4}$ and so  $ H^{\perp} = H^{\perp}_{1} \times H^{\perp}_{2} \times H^{\perp}_{3} \times H^{\perp}_{4}$  by lemma 4.2.8  of \cite{R}. Consider
 $ H_{1} < \Lambda_{1} = \mathbb{R} \times \mathbb{Z}^{n-1}$, thus by example \ref{ex3.11}, $ \frac{\Lambda_{1}}{H_{1}}$  is uncountable and consequently $ \frac{\Lambda}{H}$  is uncountable. A similar  discussion shows that $ \frac{\Gamma}{H^{\perp}} $ is uncountable, as well. Therefore, for every closed subgroup $H$ which satisfies $ \Gamma^{\perp} \leq H \leq \Lambda$, we obtain atleast one of   $ \frac{\Lambda}{H}$ or $ \frac{\Gamma}{H^{\perp}} $ are uncountable.

Now, we investigate some sufficient condition for the existence of subgroup $H$  satisfies $(\ref{existence H})$ so that $\frac{\Lambda}{H}$ and  $\frac{\Gamma}{H^{\perp}}$  are finite or countable.


\begin{theorem}
Let $G$ be an LCA group. Also, let  $\Lambda$ and $\Gamma$ be  subgroups of $G$ and $\widehat{G}$, respectively so that $\Gamma^{\perp} \leq \Lambda$. Then the following assertions hold;
\begin{itemize}
\item[(i)]
If $\Lambda$ is discrete and $\Gamma$ is  a  closed  co-compact  subgroup, then  for every  subgroup $H$
where  $\Gamma^{\perp} \leq H \leq \Lambda$,
the quotient groups  $\frac{\Lambda}{H}$  and    $\frac{\Gamma}{H^{\perp}}$  are finite.
\item[(ii)]
If   $\Lambda$ and $\Gamma$ are open subgroups, then  there exists a closed  subgroup $H$ satisfies (\ref{existence H})  so that  either $ \frac{\Lambda}{H}$ or $ \frac{\Gamma}{H^{\perp}}$ is countable.
\item[(iii)]
If $G$ is totally-disconnected and $\Lambda$, $\Gamma$ are open subgroups, then  there exists a compact subgroup $H$ satisfies (\ref{existence H})  so that both  $ \frac{\Lambda}{H}$ and $ \frac{\Gamma}{H^{\perp}}$ are countable.
\end{itemize}
\end{theorem}
\begin{proof}
$(i)$. Consider a  subgroup $H$ of $G$ such that $ \Gamma^{\perp} \leq H \leq \Lambda$. Then, $H$ is  a discrete subgroup, moreover the assumption assure  that  $\Lambda$,  $\Gamma$ and $H$ are  uniform lattices. On the other hand by Proposition 4.2.24 of \cite{R} we have $\widehat{(\frac{H^{\perp}}{\Lambda^{\perp}})} \cong \frac{\Lambda}{H}$ and so the duality relationships $ \widehat{(\frac{G}{\Lambda})} \cong \Lambda^{\perp} $ and $ \widehat{\Lambda} \cong \frac{\widehat{G}}{\Lambda^{\perp}}$  imply that $ \frac{\Lambda}{H}$ is both compact and discrete. Hence $\frac{\Lambda}{H}$ is finite, similarly the quotient group $\frac{\Gamma}{H^{\perp }}$ is finite.

$(ii)$. Since $ \Gamma$ is an open subgroup of $ \widehat{G}$, the duality relation $ \Gamma^{\perp} \cong \widehat{(\frac{\widehat{G}}{\Gamma})}$ implies that $ \Gamma^{\perp}$ is compact. So, by proposition 3.1.5 of \cite{R}, there exists a unit neighborhood $V$ of $e$ (the identity of $G$) such that $ \Gamma^{\perp} +V \subseteq \Lambda$. If we take $H:=\Gamma^{\perp} +\langle V\rangle$, then $H$ satisfies (\ref{existence H}) and    is open in $\Lambda$. Thus, $ \frac{\Lambda}{H}$ is discrete and countable but  $ \frac{\Gamma}{H^{\perp}}$ is not necessarily countable.
 On the other hand,  the structure of $V$ assure that there exists a compact subgroup $N$ of $G$ such that $ N \subseteq V$, \cite{E-K-A}. Set $H:= \Gamma^{\perp} +  N $. Then $H$ is a  compact subgroup of $G$ and  $ \Gamma^{\perp} \leq  H \leq  \Lambda$, i.e., $H$ satisfies (\ref{existence H}). Moreover,
\begin{equation*}
\dfrac{\widehat{G}}{H^{\perp}} = \dfrac{\widehat{G}}{\Gamma \cap N^{\perp}} \cong \widehat{(\Gamma^{\perp} + N)}
\end{equation*}
is discrete and so countable. Therefore, $ \frac{\Gamma}{H^{\perp}}$ is also countable, but not necessarily $ \frac{\Lambda}{H}$.

$(iii)$. Using the proof of  $(ii)$  there exists a unit neighborhood $V$ of $e$ so that $ \Gamma^{\perp} +V \subseteq \Lambda$. The assumption that $G$  is totally-disconnected implies that
  there exists  an open compact subgroup  $K$  so that $K \subseteq V$,  by Theorem 7.7 of \cite{E-K-A}. Put $H:= \Gamma^{\perp} + K$. Then $H$ is a compact open  subgroup of $G$ and  $ \Gamma^{\perp} \leq  H  \leq \Lambda$. Moreover, $H$ is open in $\Lambda$. Thus, $ \frac{\Lambda}{H}$ is   countable. Also, an analogous discussion shows that $ \frac{\Gamma}{H^{\perp}}$ is  countable as well. This completes the proof.
\end{proof}


\section{Fibrization method}
The fiberization technique is closely related to Zak transform methods in Gabor analysis. Let $H$ be a closed and co-compact subgroup of $G$ and $\Omega \subset \widehat{G}$ be a Borel section of $H^{\perp}$ in $\widehat{G}$,  we consider the fiberization mapping $\mathcal{T} : L^{2}(G) \longrightarrow L^{2}(\Omega, l^{2}(H^{\perp}))$, introduced in \cite{B-R} by
\begin{equation*}
\mathcal{T} f(\omega) = \lbrace \widehat{f}(\omega + \alpha) \rbrace_{\alpha \in H^{\perp}}, \ \ \ \  (\omega \in \Omega).
\end{equation*}
The fiberization is a isometric isomorphic operation as shown in \cite{B-R}. Furthermore,
 the frame property of translation-invariant and Gabor system can be characterized in terms of fibers \cite{J-L}.
\begin{theorem}\label{thm}
Let $A$ and $B$ be two positive constants and let $H\subseteq G$ be a closed, co-compact subgroup, and let $ \lbrace g_{p} \rbrace_{p \in P} \subseteq L^{2}(G)$, where $(P, \mu_{P})$ is an admissible measure space. Then the following assertions are equivalent.

(i) The family $\lbrace T_{h} g_{p} \rbrace_{h \in H, p \in P}$ is a frame for $L^{2}(G)$ with bounds $A$ and $B$.

(ii) For almost every $\omega \in \Omega$, the family $\lbrace \mathcal{T} g_{p}(\omega) \rbrace_{p \in P}$ is a frame for $ l^{2}(H^{\perp})$ with bounds $A$ and $B$, where $ \Omega$ is a Borel section of $H^{\perp}$ in $ \widehat{G}$.
\end{theorem}

The next result shows that the frame property of a Gabor system in $ L^{2}(G)$ under certain assumptions is equivalent with the frame property of a family of associated Zak transforms in $ l^{2}(\widehat{H^{\perp}})$.
\begin{theorem}
Let $g \in L^{2}(G)$,  $\Lambda$ and $\Gamma$ be closed subgroups of $G$ and $\widehat{G}$ respectively and let  $H$ be a closed, co-compact subgroup of $G$ satisfies \eqref{existence H}. Then there exists a sequence $ \lbrace  g_{ku} \rbrace_{ k \in \frac{\Lambda}{H} , u \in \frac{\Gamma}{\Lambda^{\perp}}}$ in $ L^{2}(G)$ such that following assertions are equivalent.
\begin{itemize}
\item[(i)] $ \lbrace E_{\gamma} T_{\lambda} g \rbrace_{\lambda \in \Lambda , \gamma \in \Gamma}$ is a frame for $L^{2}(G)$ with bounds $A$ and $B$.
\item[(ii)] $ \lbrace Z_{H^{\perp}} \widehat{ g_{ku}} ( \omega , .) \rbrace_{ k \in \frac{\Lambda}{H} , u \in \frac{\Gamma}{\Lambda^{\perp}}} $ is a frame for $l^{2}(\widehat{H^{\perp}})$ with bounds $A$ and $B$,  for a.e.  $ \omega \in \Omega$, where $ \Omega$ is a Borel section of $H^{\perp}$ in $ \widehat{G}$.
\end{itemize}
\end{theorem}
\begin{proof}
Since $ H \leq \Lambda$, so  every $ \lambda \in \Lambda$ can be written uniquely as  $ \lambda = t +k$ where $t \in H$ and $ k \in \frac{\Lambda}{H}$. Also, $ \Lambda^{\perp} \leq \Gamma$    implies that every $ \gamma \in \Gamma$ has a unique form  such as $ \gamma = \mu + u$ where $ \mu \in \Lambda^{\perp}$ and $ u \in \frac{\Gamma}{\Lambda^{\perp}}$. Thus,
\begin{eqnarray*}
 \lbrace T_{\lambda} E_{\gamma} g \rbrace_{\lambda \in \Lambda , \gamma \in \Gamma}
= \lbrace T_{t} E_{\mu} g_{ku} \rbrace_{t \in H , k \in \frac{\Lambda}{H} , \mu \in \Lambda^{\perp} , u \in \frac{\Gamma}{\Lambda^{\perp}}}
\end{eqnarray*}
 where $\mathcal{G}:=\lbrace g_{ku} \rbrace = \lbrace  T_{k} E_{u} g \rbrace_{ k \in \frac{\Lambda}{H} , u \in \frac{\Gamma}{\Lambda^{\perp}}}$.
Therefore,  applying the fiberization method along with Theorem \ref{thm}  for co-compact subgroup $H$ of $G$, the   system $ \lbrace T_{t} E_{\mu} \mathcal{G} \rbrace_{ t \in H, \mu \in \Lambda^{\perp}}$ (or equivalently $\lbrace E_{\gamma} T_{\lambda} g \rbrace_{\lambda \in \Lambda , \gamma \in \Gamma}$) is a frame for $L^{2}(G)$ if and only if $ \lbrace \mathcal{T} E_{\mu} \mathcal{G} (\omega) \rbrace_{ \mu \in \Lambda^{\perp}}$ is a frame in $ l^{2}(H^{\perp})$, for  a.e $ \omega \in \Omega$ where $ \Omega$ is a Borel section of $H^{\perp}$ in $ \widehat{G}$. On the other hand,  we obtain
\begin{eqnarray*}
\lbrace \mathcal{T} E_{\mu} \mathcal{G}(\omega) \rbrace_{ \mu \in \Lambda^{\perp}}
&=& \lbrace \mathcal{T} E_{\mu} g_{ku}(\omega) \rbrace_{ \mu \in \Lambda^{\perp} , k \in \frac{\Lambda}{H} ,  u \in \frac{\Gamma}{\Lambda^{\perp}}}\\
&=& \lbrace  \lbrace  \widehat{E_{\mu} g_{ku}} ( \omega + \alpha) \rbrace_{ \alpha \in H^{\perp}} \rbrace_{ \mu \in \Lambda^{\perp} , k \in \frac{\Lambda}{H} ,  u \in \frac{\Gamma}{\Lambda^{\perp}}} \\
&=& \lbrace \lbrace T_{\mu} \widehat{g_{ku}} ( \omega + \alpha ) \rbrace_{ \alpha \in H^{\perp}} \rbrace_{\mu \in \Lambda^{\perp} , k \in \frac{\Lambda}{H} ,  u \in \frac{\Gamma}{\Lambda^{\perp}}}\\
&=& \lbrace  \lbrace  \widehat{g_{ku}}(\omega + \alpha ) \rbrace_{ \alpha \in H^{\perp}} \rbrace_{ k \in \frac{\Lambda}{H} ,  u \in \frac{\Gamma}{\Lambda^{\perp}}},
\end{eqnarray*}
where the last equality is because of  the assumption $(\ref{existence H})$.
 Consider $ \psi_{ k,u} (\omega) := \lbrace \widehat{g_{ku}} ( \omega + \alpha ) \rbrace_{ \alpha \in H^{\perp}}$ for all $ k \in \frac{\Lambda}{H} ,  u \in \frac{\Gamma}{\Lambda^{\perp}} $ and a.e $ \omega \in \Omega $. Then  the Fourier inversion transform of $ \psi_{ k,u}(\omega) \in l^{2}(H^{\perp})$ is as follows
\begin{eqnarray*}
\mathcal{F}^{-1} (\psi_{k,u}(\omega)) (\xi) &=& \sum_{ \alpha \in H^{\perp}} \widehat{g_{ku}} ( \omega + \alpha) \alpha (\xi)\\
&=&Z_{H^{\perp}} \widehat{g_{ku}} ( \omega  , \xi)
\end{eqnarray*}
for all $k \in \frac{\Lambda}{H}$ , $ u \in \frac{\Gamma}{\Lambda^{\perp}}$  and  a.e., $\xi \in \widehat{H^{\perp}}$  and $\omega \in \Omega$.
Hence, the assertion $(i)$ is equivalent to the system   $\lbrace Z_{H^{\perp}} \widehat{g_{ku}} ( \omega ,.) \rbrace_{ k \in \frac{\Lambda}{H} ,  u \in \frac{\Gamma}{\Lambda^{\perp}}}$ being a frame for $l^{2}(\widehat{H^{\perp}})$, a.e  $ \omega \in \Omega$, as required.
\end{proof}
In the next corollary, we exploit some connections to the results obtained in \cite{J-L}.
\begin{corollary}
Let $g \in L^{2}(G)$,  $\Lambda$ be a  closed co-compact subgroup  of $G$   and $\Gamma$ be a  closed subgroup  of  $\widehat{G}$ so that $\Gamma^{\perp}\leq \Lambda$. Then following assertions are equivalent.
\begin{itemize}
\item[(i)]$ \lbrace E_{\gamma} T_{\lambda} g \rbrace_{\lambda \in \Lambda , \gamma \in \Gamma}$ is a frame for $L^{2}(G)$ with bounds $A$ and $B$.
\item[(ii)]$\{\widehat{g}(\alpha+\gamma)\}_{\gamma\in \Gamma}$ is a frame for $l^{2}(\Lambda^{\perp})$  with bounds $A$ and $B$, for a.e. $\alpha\in \mathcal{A}$, where $\mathcal{A}$ is a Borel section of $\Lambda^{\perp}$ in $\widehat{G}$.
\item[(iii)] $A\leq \int_{\mathcal{K}}  \left\vert Z_{\Lambda^{\perp}}\widehat{g}(\alpha+ k,x)\right\vert^{2} d\mu_{\mathcal{K}}{(k)}   \leq B$, for a.e. $\alpha \in \mathcal{A}$ and $x\in \widehat{\Lambda^{\perp}}$, where $\mathcal{A}$  is a Borel section of $\Lambda^{\perp}$ in $\widehat{G}$, $\mathcal{K}\subset \Gamma$  is a Borel section of $\Lambda^{\perp}$ in $\Gamma$.
\end{itemize}
\end{corollary}
\begin{proof}
$(i)\Leftrightarrow(ii)$; We note that $ (\Gamma, \Sigma_{\Gamma}, \mu_{\Gamma})$ is an admissible measure space, since $\Gamma$ is a closed subgroup of  $\widehat{G}$.  Hence, by Proposition 4.5 in \cite{J-L}, $(i)$ is equivalent to $ \lbrace \lbrace \widehat{g} (\alpha +\gamma+ y ) \rbrace_{y \in \Lambda^{\perp}} \rbrace_{\gamma \in \Gamma}= \{\widehat{g}(\alpha+\gamma)\}_{\gamma\in \Gamma}$ being a frame for $ l^{2}(\Lambda^{\perp})$ with bounds $A$ and $B$, for a.e. $ \alpha \in \mathcal{A} $, where $ \mathcal{A}$ is a Borel section of $ \Lambda^{\perp}$ in $\widehat{G}$.

$(i)\Leftrightarrow(iii)$ Applying   the fact that  $\Lambda$ is co-compact  and
  $\Lambda^{\perp}\cap \Gamma=\Lambda^{\perp}$, it is sufficient to take appropriate Haar measures on $\Gamma$, $\frac{\Gamma}{\Lambda^{\perp}}$  and put a measure on $\mathcal{K}$ isometric to $\mu_{\Gamma/\Lambda^{\perp}}$ in the sense of $(\ref{measure})$. Then the desired result is  obtained by  Theorem 4.6 in \cite{J-L}.
\end{proof}


\end{document}